\documentclass{IET-Conf-Paper}
\usepackage[utf8]{inputenc}
\usepackage{balance}
\usepackage{hyperref}
\DeclareMathSizes{10}{9}{7}{5}
\usepackage{amssymb}
\usepackage{nomencl}
\usepackage{etoolbox}
\usepackage{graphics}
\usepackage{url}
\usepackage{tabularx}
\allowdisplaybreaks
\usepackage{epstopdf}
\usepackage{color}
\usepackage{xcolor,colortbl}
\usepackage{tabularx,colortbl}
\usepackage{balance}
\usepackage{amssymb}
\usepackage{graphicx}
\usepackage{adjustbox}
\usepackage{natbib}
\usepackage{dirtytalk}
\usepackage{color}
\usepackage{xcolor,colortbl}
\usepackage{tabularx,colortbl}
\usepackage{pgfplots}
\usepackage{pgfplotstable}
\usetikzlibrary{decorations.pathmorphing}

\pgfplotsset{compat=newest}
\usepackage{caption}
\usepackage{subcaption}
\usepackage{subfiles}
\usepackage{standalone}
\usepackage{tikz}
\usepackage{mathtools,nccmath}
\newcolumntype{P}[1]{>{\centering\arraybackslash}p{#1}}
\begin{document}

\title{Towards a Sustainable Microgrid on Alderney Island Using a Python-based Energy Planning Tool}

\author{Shahab Dehghan\ad{1}\corr, Agnes M Nakiganda\ad{1}, James Lancaster\ad{2}, Petros Aristidou\ad{1,3}}

\address{\add{1}{School of Electronic and Electrical Engineering, 
 University of Leeds, Leeds, UK}
\add{2}{Alderney Electricity Ltd, Alderney, Channel Islands}
\add{3}{Department of Electrical Engineering, Computer Engineering \& Informatics,\\Cyprus University of Technology, Cyprus}
\email{s.dehghan@leeds.ac.uk}}
\vspace{-2cm}
\vspace{-0.5cm}
\begin{abstract}
\vspace{-0.25cm}
In remote or islanded communities, the use of microgrids (MGs) is necessary to ensure electrification and resilience of supply. However, even in small-scale systems, it is computationally and mathematically challenging to design low-cost, optimal, sustainable solutions taking into consideration all the uncertainties of load demands and power generations from renewable energy sources (RESs). This paper uses the open-source Python-based Energy Planning (PyEPLAN) tool, developed for the design of sustainable MGs in remote areas, on the Alderney island, the 3$^{rd}$ largest of the Channel Islands with a population of about 2000 people. A two-stage stochastic model is used to optimally invest in battery storage, solar power, and wind power units. Moreover, the AC power flow equations are modelled by a linearised version of the DistFlow model in PyEPLAN, where the investment variables are {\it here-and-now} decisions and not a function of uncertain parameters while the operation variables are {\it wait-and-see} decisions and a function of uncertain parameters. The $k$-means clustering technique is used to generate a set of best (risk-seeker), nominal (risk-neutral), and worst (risk-averse) scenarios capturing the uncertainty spectrum using the yearly historical patterns of load demands and solar/wind power generations. The proposed investment planning tool is a mixed-integer linear programming (MILP) model and is coded with Pyomo in PyEPLAN.
\end{abstract}

\maketitle

\thanksto{\noindent This work is supported by the UK Engineering and Physical Sciences Research Council (EPSRC) under Grant EP/R030243/1.}

\makenomenclature
\renewcommand\nomgroup[1]{%
  \item[\bfseries
  \ifstrequal{#1}{A}{Indices}{%
  \ifstrequal{#1}{B}{Parameters}{%
  \ifstrequal{#1}{C}{Sets}{%
  \ifstrequal{#1}{D}{Variables}{}}}}%
]}

\nomenclature[A, 01]{$n$}{Index of nodes where $n'$ and $n''$ stand for nodes before and after node $n$, respectively.}
\nomenclature[A, 02]{$d$}{Index of load demands.}
\nomenclature[A, 03]{$g$}{Index of generation units.}
\nomenclature[A, 04]{$o$}{Index of representative days (scenarios).}
\nomenclature[A, 05]{$t$}{Index of time periods.}

\nomenclature[B, 01]{$e_{b}^{\textrm{ini}}$}{Initial stored energy of battery unit $b$ (kW).}
\nomenclature[B, 02]{$e_{b}^{\max}$}{Maximum stored energy of battery unit $b$ (kW).}
\nomenclature[B, 03]{$e_{b}^{\min}$}{Minimum stored energy of battery unit $b$ (kW).}
\nomenclature[B, 04]{$pc_{d}$}{Penalty cost of load demand curtailment (\$/kWh).}
\nomenclature[B, 05]{$pc_{r}$}{Penalty cost of RES power generation curtailment (\$/kWh).}
\nomenclature[B, 06]{$f_{d}$}{Power factor of load demand $d$.}
\nomenclature[B, 07]{${ic}_{b}$}{Annualised investment cost of battery unit $b$ (\$).}
\nomenclature[B, 08]{${ic}_g$}{Annualised investment cost of generation unit $g$ (\$).}
\nomenclature[B, 09]{$mc_{g}$}{Marginal cost of generation unit $g$ (\$/kWh).}
\nomenclature[B, 10]{$p_{b}^{\max,\textrm{c/d}}$}{Maximum charging/discharging power of battery unit $b$ (kW).}
\nomenclature[B, 11]{$\bar p_{dto}$}{Load demand $d$ at hour $t$ in representative day $o$ (\$/kWh).}
\nomenclature[B, 12]{$p_{n'n}^{\max}$}{Maximum active power flow from node $n'$ to node $n$ (kW).}
\nomenclature[B, 13]{$p_{g}^{\max}$}{Maximum active power of generation unit $g$ (kW).}
\nomenclature[B, 14]{$\overline p_{gto}^{\max}$}{Maximum power generation of generation unit $g$ at hour $t$ in representative day $o$ (kW).}
\nomenclature[B, 15]{$q_{b}^{\max}$}{Maximum reactive power of battery unit $b$ (kVAr).}
\nomenclature[B, 16]{$q_{b}^{\min}$}{Minimum reactive power of battery unit $b$ (kVAr).}
\nomenclature[B, 17]{$q_{n'n}^{\max}$}{Maximum reactive power flow from node $n'$ to node $n$ (kVAr).}
\nomenclature[B, 18]{$q_{g}^{\max}$}{Maximum reactive power of generation unit $g$ (kVAr).}
\nomenclature[B, 19]{$q_{g}^{\min}$}{Minimum reactive power of generation unit $g$ (kVAr).}
\nomenclature[B, 20]{$r_{n'n}$}{Resistance of the line connecting nodes $(n',n)$ (ohm).}
\nomenclature[B, 21]{$v^{\max}$}{Maximum permitted voltage magnitude (V).}
\nomenclature[B, 22]{$v^{\min}$}{Minimum permitted voltage magnitude (V).}
\nomenclature[B, 23]{$x_{n'n}$}{Reactance of the line connecting nodes $(n',n)$ (ohm).}
\nomenclature[B, 24]{$\eta_b^{\textrm{c/d}}$}{Reactance of the line connecting nodes $(n',n)$ (ohm).}

\nomenclature[C, 01]{$\Omega^{B}$}{Set of battery units where $\Omega^{B_n}$ indicates set of battery units connected to node $n$.}
\nomenclature[C, 02]{$\Omega^{D}$}{Set of load demands where $\Omega^{D_n}$ indicates set of load demands connected to node $n$.}
\nomenclature[C, 03]{$\Omega^{L}$}{Set of distribution lines connecting nodes.}
\nomenclature[C, 04]{$\Omega^{M}$}{Set of micro-turbine/diesel units where $\Omega^{M_n}$ indicates set of micro-turbine/diesel generators connected to node $n$.}
\nomenclature[C, 01]{$\Omega^{N}$}{Set of nodes where $\Omega^{N_n}$ indicates set of nodes after and connected to node $n$.}
\nomenclature[C, 05]{$\Omega^{R}$}{Set of RES units where $\Omega^{R_n}$ indicates set of RES units connected to node $n$.}
\nomenclature[C, 06]{$\Omega^{T}$}{Set of hours.}

\nomenclature[D, 01]{$p_{bto}^{\textrm{c/d}}$}{Active charging/discharging power of battery unit $b$ at hour $t$ in representative day $o$ (kW).}
\nomenclature[D, 02]{$p_{n'nto}$}{Active power flow from node $n'$ to node $n$ at hour $t$ in representative day $o$ (kW).} 
\nomenclature[D, 03]{$p_{gto}$}{Active power generation of generation unit $g$ at hour $t$ in representative day $o$ (kW).}
\nomenclature[D, 04]{$q_{bto}$}{Reactive power of battery unit $b$ at hour $t$ in representative day $o$ (kW).}
\nomenclature[D, 05]{$q_{n'nto}$}{Reactive power flow from node $n'$ to node $n$ at hour $t$ in representative day $o$ (kVAr).}
\nomenclature[D, 06]{$q_{gto}$}{Reactive power generation of generator $g$ at hour $t$ in representative day $o$ (kVAr).}
\nomenclature[D, 08]{$v_{nto}$}{Voltage magnitude of node $n$ at hour $t$ in representative day $o$ (V).}
\nomenclature[D, 10]{$y_{dto}$}{Curtailment status of load demand $d$ at hour $t$ in representative day $o$ (i.e., 1/0: curtailed/not-curtailed).} 
\nomenclature[D, 10]{$z_{b}$}{Investment status of battery unit $b$ (i.e., 1/0: built/non-built).} 
\nomenclature[D, 11]{$z_{g}$}{Investment status of RES unit $g$ (i.e., 1/0: built/non-built).}
\vspace{-1.5cm}
\printnomenclature
\section{Introduction}\label{introduction}
Alderney island with an area of 3 square miles runs a closed complex energy system that entirely relies on imported fuel oils for electricity, heating, and transportation. Major economic activities on the island include e-trade, ecotourism, small businesses, health care services. The only energy supplier on the island is Alderney Electricity Limited (AEL)~\cite{AEL_Website}, providing for both electric and heating loads. AEL is responsible for the importation and distribution of different fuels, including kerosene and transport fuels, as well as the generation and distribution of electricity. The company manages both the $11$ kV primary distribution network, consisting of $21$ substations, as well as the $415$ V secondary distribution network. AEL starts with the higher voltage to account for cable losses ensuring the voltage is still in spec. by the time it gets where it is going. Networks mainly comprise underground cables, there are a small number of overheads which are being progressively replaced. Electric power on Alderney island is centrally generated by $8\times450$ kVA diesel generators and supplied through an extensive network consisting of underground cables. Hence, the main aim of this paper is to create a sustainable microgrid (MG) on Alderney island, which obviates the reliance of AEL on only fossil fuels.

\vspace{-0.25cm}
\subsection{Literature Review and Contributions}\label{Literature-Review}
\vspace{-0.25cm}
MG is a low-voltage electrical network, including diverse controllable and uncontrollable producers, consumers, and prosumers, that can be operated autonomously. The concept of MG has been initially introduced in the seminal reference~\cite{Microgrids-distributed-power-generation} to cope with the main challenges in integrating distributed energy resources into low-voltage electric networks. Most of MGs in remote areas (like Alderney island) have been operated by fossil fuel-based generation technologies with competitive costs as compared to sustainable generation technologies. However, increasing concerns related to global climate change as well as advances in sustainable generation technologies have made renewable energy sources (RESs) a priority in MGs during the last decade~\cite{Trends-in-Microgrid-Control}. Since RES power generation (e.g., solar and wind power) is inherently subject to uncertainty and volatility, ignoring them may result in infeasible investment and operation plans. Therefore, it is of utmost importance to use practical investment and operation planning tools presenting feasible solutions under different uncertainties. 

Previously, stochastic optimisation (SO)~\cite{Microgrid-Planning-Under-Uncertainty,Provisional-Microgrid-Planning} and robust optimisation (RO)~\cite{Stochastic-Capacity-Expansion-Planning-of-Remote-Microgrids-With-Wind-Farms-and-Energy-Storage,Integrated-Microgrid-Expansion-Planning-in-Electricity-Market-With-Uncertainty,amjady2017adaptive} have been introduced in the literature to cope with different uncertain parameters in distribution networks and MGs. RO provides an investment plan, which is optimal under the worst-case scenario of uncertain parameters, while SO provides an investment plan, which is optimal on average for all scenarios characterising uncertain parameters. It is noteworthy to mention that the optimal solutions of RO-based investment planning models may be conservative than the optimal solutions of SO-based ones in MGs with sufficient historical data. Accordingly, a Python-based Energy PLANning (PyEPLAN) tool is used in this paper to propose a sustainable MG strategy on Alderney island based on a two-stage SO-based model. In the proposed approach, investment variables are {\it here-and-now} decisions and not a function of uncertain parameters, while operation variables are {\it wait-and-see} decisions and a function of uncertain parameters. In summary, the main contributions of this paper are as follows: (i) A two-stage stochastic mixed-integer linear programming (MILP) model is introduced in this paper to optimally invest in battery, solar, and wind units on Alderney island under the uncertainty of load demands and RES power generations; (ii) A practical MG test system is presented for future investment and operation planning studies based on the network data of the AEL MG.  
\vspace{-0.25cm}
\subsection{Paper Organisation}\label{Literature-Review}
\vspace{-0.25cm}
The rest of this paper is organised as follows. In Section~\ref{section2}, an overview about PyEPLAN and its clustering, investment planning, and operation planning modules are presented. In Section~\ref{section3}, the proposed two-stage stochastic MG investment planning (SMIP) model as an MILP optimisation problem is introduced. In Section~\ref{section4}, the proposed SMIP model is tested on the AEL MG under different conditions. Finally, Section~\ref{section5} concludes the paper.    

\section{Brief Review of PyEPLAN}\label{section2}
\begin{figure}[!t]
\centering
\includegraphics[width=0.5\textwidth]{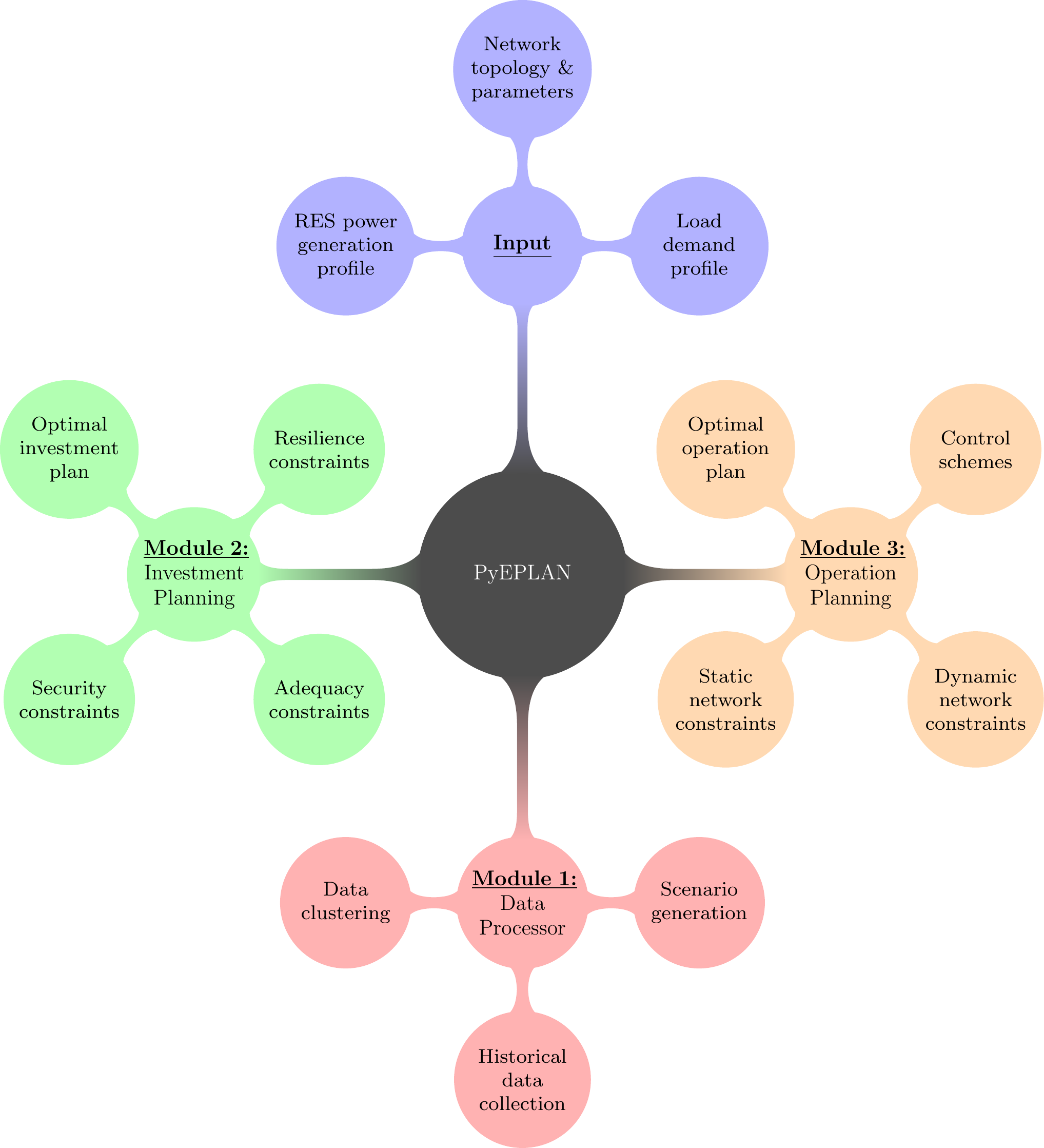}  
\caption{Overall PyEPLAN architecture.}
\label{Fig: PyEPLAN}
\end{figure}

The planning tool used in this work, PyEPLAN, has three different modules, including data processing, investment planning, and operation planning in MGs, as depicted in Fig.~\ref{Fig: PyEPLAN}. In this paper, only the investment planning module is used to plan a sustainable MG on Alderney island. Internally, PyEPLAN uses the open-source Python-based optimisation modelling (Pyomo)~\cite{hart2017pyomo} language to formulate, solve, and analyze the optimisation problems for investment and operation planning. 
Both investment and operation planning modules in PyEPLAN are developed based on a concrete~\cite{hart2017pyomo} model of Pyomo that can be initialised by means of comma-separated values (CSV) files, including input data sets (i.e., different characteristics of various components in MGs). 

\vspace{-0.25cm}
\subsection{Investment Planning Module}\label{SMIP} 
\vspace{-0.25cm}
The objective of the MILP is to minimize both investment and operation costs during a {\it long-term} planning horizon (i.e., from one year to several years) under both investment and operation related techno-economic constraints. As input, the module needs network characteristics (i.e., candidate/existing generation technologies, candidate/existing lines), as well as long-term estimated/forecasted load demands and RES power generations to obtain the optimal solution. Accordingly, the data processor, as discussed in the next subsection, is considered in PyEPLAN to provide the input data needed for the investment planning module. 
\vspace{-0.25cm}
\subsection{Data Processor}\label{kmeans}
\vspace{-0.25cm}
In the investment planning module, it is assumed that the pattern of load demands (obtained by dividing the hourly load demands of each year by its peak), as well as the pattern of RES power generations (obtained by dividing the hourly power generations of each RES by its capacity) remain unchanged during a one-year period~\cite{dehghan2015reliability}. However, the SMIP model needs a sufficient number of scenarios to characterise the uncertain load demand as well as the uncertain RES power generation during a one-year period. Therefore, the $k$-means clustering technique, as presented in~\cite{dehghan2019robust}, is used to obtain representative days from daily load demand profiles and RES power generation during a year. Then, the SMIP model incorporates the best (risk-seeker), nominal (risk-neutral), and worst (risk-averse) representative days~\cite{dehghan2015reliability}. 

\section{Stochastic MG Planning Model}\label{section3}
In this section, the mathematical formulation of the proposed SMIP model is briefly reviewed within a single-year planning horizon under different representative days (scenarios) for load demands and RES power generations as given below: 
\begin{subequations}
\begin{align}
\label{1a}
\begin{split}
\min \mathbf{\Psi}^{\textrm{inv}} + \mathbf{\Psi}^{\textrm{opr}}
\end{split}
\end{align}
\text{s.t.}
\begin{align}
\begin{split}
\label{1b}
&\mathbf{\Psi}^{\textrm{inv}} = \sum_{b\in \Omega^{B} }\left(ic_{b}\cdot z_{b}\right) + \sum_{g\in \Omega^{R} }\left(ic_{g}\cdot z_{g}\right) 
\end{split}\\
\begin{split}
\label{1c}
&\mathbf{\Psi}^{\textrm{{opr}}}=\sum_{o\in\Omega^{O}}\sum_{t\in\Omega^{T}}\sum_{g\in\{\Omega^{M},\Omega^{R}\}}\left(\tau_{o}\cdot mc_{g}\cdot p_{gto}\right)+\\
&\sum_{o\in\Omega^{O}}\sum_{t\in\Omega^{T}}\sum_{s\in\Omega^{S}}\left(\tau_{o}\cdot pc_{d}\cdot\overline{p}_{dto}\cdot \left(1-y_{dot}\right)\right)+\\
&\sum_{o\in\Omega^{O}}\sum_{t\in\Omega^{T}}\sum_{g\in\Omega^{R}}\left(\tau_{o}\cdot pc_{r}\cdot\left(\overline{p}_{gto}^{\max}-p_{gto}\right)\right)
\end{split}\\
\begin{split}
\label{1d}
&p_{n^{\prime}nto}+\sum_{g\in\{\Omega^{M_{n}},\Omega^{R_{n}}\}}p_{gto}+\sum_{b\in\Omega^{B_{n}}}\left(p_{bto}^{\textrm{d}}-p_{bto}^{\textrm{c}}\right)=\\
&\sum_{n^{\prime\prime}\in\Omega^{N_{n}}}p_{nn^{\prime\prime}to}+\sum_{d\in\Omega^{D_{n}}}\left(\bar{p}_{dto} \cdot y_{dto}\right)\,\,\, n\in\Omega^{N},t\in\Omega^{T},o\in\Omega^{O}
\end{split}\\
\begin{split}\label{1e}
&q_{n^{\prime}nto}+\sum_{g\in\Omega^{M_{n}}}q_{gto}+\sum_{b\in\Omega^{B_{n}}}q_{bto}=\sum_{n^{\prime\prime}\in\Omega^{N_{n}}}q_{nn^{\prime\prime}t}+\\
&\sum_{d\in\Omega^{D_{n}}}\tan\left(\arccos\left(f_{d}\right)\right)\cdot\left(\bar{p}_{dto}\cdot y_{dto}\right)\,\,\, n\in\Omega^{N},t\in\Omega^{T},o\in\Omega^{O}
\end{split}\\
\begin{split}
\label{1f}
&\left(r_{n^{\prime}n}\cdot p_{n^{\prime}nto}+x_{n^{\prime}n}\cdot q_{n^{\prime}nto}\right)=\\
&v_{n^{\prime}to}-v_{nto} \quad n\in\Omega^{N},t\in\Omega^{T},o\in\Omega^{O}
\end{split}\\
\begin{split}\label{1g}
&-p_{nn^{\prime\prime}}^{\max}\le p_{nn^{\prime\prime}to}\le p_{nn^{\prime\prime}}^{\max} \,\left(n,n^{\prime\prime}\right)\in\Omega^{L},t\in\Omega^{T},o\in\Omega^{O}
\end{split}\\
\begin{split}
\label{1h}
&-q_{nn^{\prime\prime}}^{\max}\le q_{nn^{\prime\prime}to}\le q_{nn^{\prime\prime}}^{\max}\,\, \left(n,n^{\prime\prime}\right)\in\Omega^{L},t\in\Omega^{T},o\in\Omega^{O}
\end{split}\\
\begin{split}
\label{1i}
&0\leq p_{gto}\le p_{g}^{\max} \quad g\in\Omega^{M},t\in\Omega^{T},o\in\Omega^{O}
\end{split}\\
\begin{split}
\label{1j}
&q_{g}^{\min} \leq q_{gto}\le q_{g}^{\max} \quad g\in\Omega^{M},t\in\Omega^{T},o\in\Omega^{O}
\end{split}\\
\begin{split}
\label{1k}
&0\leq p_{gt}\le\bar{p}_{gto}^{\max}\cdot z_{g}\quad g\in\Omega^{R},t\in\Omega^{T},o\in\Omega^{O}
\end{split}\\
\begin{split} \label{1l}
& q_{g}^{\min}\cdot z_{g}\leq q_{gto}\le q_{g}^{\max}\cdot z_{g}\quad g\in\Omega^{R},t\in\Omega^{T},o\in\Omega^{O}
\end{split}\\
\begin{split} \label{1m}
& {e}_b^{\min} \cdot z_b  \leq e_{bo}^{\textrm{ini}} + \sum_{\tau = 1}^{t}\left( \eta_b^{\textrm{c}} \cdot p^{{\textrm{c}}}_{b\tau o} - \dfrac{1}{\eta_b ^{\textrm{d}}}\cdot p^{{\textrm{d}}}_{b\tau o}\right) \\
&\qquad \qquad \le {e}_b^{\max}\cdot z_b \quad b\in \Omega^B,t\in\Omega^{T},o\in\Omega^{O}
\end{split}\\
\begin{split}\label{1n}
& \sum_{\tau = 1}^{T}\left( \eta_b^{{\textrm{c}}} \cdot p^{{\textrm{c}}}_{b\tau o} - \dfrac{1}{\eta ^{{\textrm{d}}}}\cdot p^{{\textrm{d}}}_{b\tau o}\right)=0 \: b\in \Omega^B,t\in\Omega^{T},o\in\Omega^{O}
\end{split}\\
\begin{split}\label{1o}
& 0\le p_{bto}^{\textrm{{c}}}\le p^{\max,\textrm{{c}}}\cdot z_b \quad b\in \Omega^B,t\in\Omega^{T},o\in\Omega^{O}
\end{split}\\
\begin{split}\label{1p}
& 0\le p_{bto}^{\textrm{{d}}}\le p^{\max,\textrm{{d}}}\cdot z_b\quad b\in \Omega^B,t\in\Omega^{T},o\in\Omega^{O}
\end{split}\\
\begin{split}\label{1q}
&v^{\min}\le v_{not}\le v^{\max}\quad n\in\Omega^{N},t\in\Omega^{T},o\in\Omega^{O}
\end{split}\\
\begin{split}\label{1r}
&v_{1to}=1\quad t\in\Omega^{T},o\in\Omega^{O}
\end{split}
\end{align}
\end{subequations}

The objective function \eqref{1a} minimises the total investment and operational costs, where $\mathbf{\Psi}^{\textrm{inv}}$ calculates the total investment costs of battery and RES units, as indicated in \eqref{1b}, and $\mathbf{\Psi}^{\textrm{opr}}$ represents the total operational costs of micro-turbine/diesel and RES units as well as curtailment costs of load demands and RES power generations, as indicated in \eqref{1c}. For simplicity, all existing and candidate technologies are considered as investment candidates, where the investment costs (resp. decision variables) of existing technologies (i.e., micro-turbine/diesel units) are set to $0$ (resp. $1$). 

PyEPLAN offers different ways to include the network constraints. In this paper, the linearised approximation of the DistFlow formulation is selected for the AC power flow equations~\cite{Network-reconfiguration-in-distribution-systems-for-loss-reduction-and-load-balancing} and the quadratic power flow limitations are linearised by means of a polygon approximation~\cite{A-Linearized-OPF-Model-With-Reactive-Power-and-Voltage-Magnitude}. Accordingly, constraints \eqref{1d} and \eqref{1e} ensure active and reactive power balance at each node of every hour of all representative days, respectively. Constraint \eqref{1f} denotes the difference of voltage magnitudes between two neighbor nodes connected. Constraints \eqref{1g} and \eqref{1h} bound the active and reactive power flows between two connected neighbor nodes, respectively. Constraints \eqref{1i} and \eqref{1j} ensure the limits on active and reactive power generation for micro-turbine/diesel units, respectively, while constraints \eqref{1k} and \eqref{1l} ensure the limits of active power generation for RES units. 

Constraint \eqref{1m} bounds the stored energy of each battery unit at every hour of all representative days. Moreover, constraint \eqref{1n} ensures the initial and final stored energy of battery units for each representative day. Constraints \eqref{1o} and \eqref{1p} bound the charging and discharging power of each battery unit at every hour in all representative days, respectively. Constraint \eqref{1q} limits the allowed variation bound of the nodal voltage magnitude. Also, constraint \eqref{1r} sets the voltage magnitude at the main AEL substation on one. The SMIP model in \eqref{1a}-\eqref{1r} is an MILP problem, which can be solved by off-the-shelf optimisation packages.
\begin{figure}[!t]
	\centering
	\scalebox{1.0}{\includegraphics[clip, trim=1.2cm 0.5cm 1.5cm 0.5cm, width=\columnwidth]{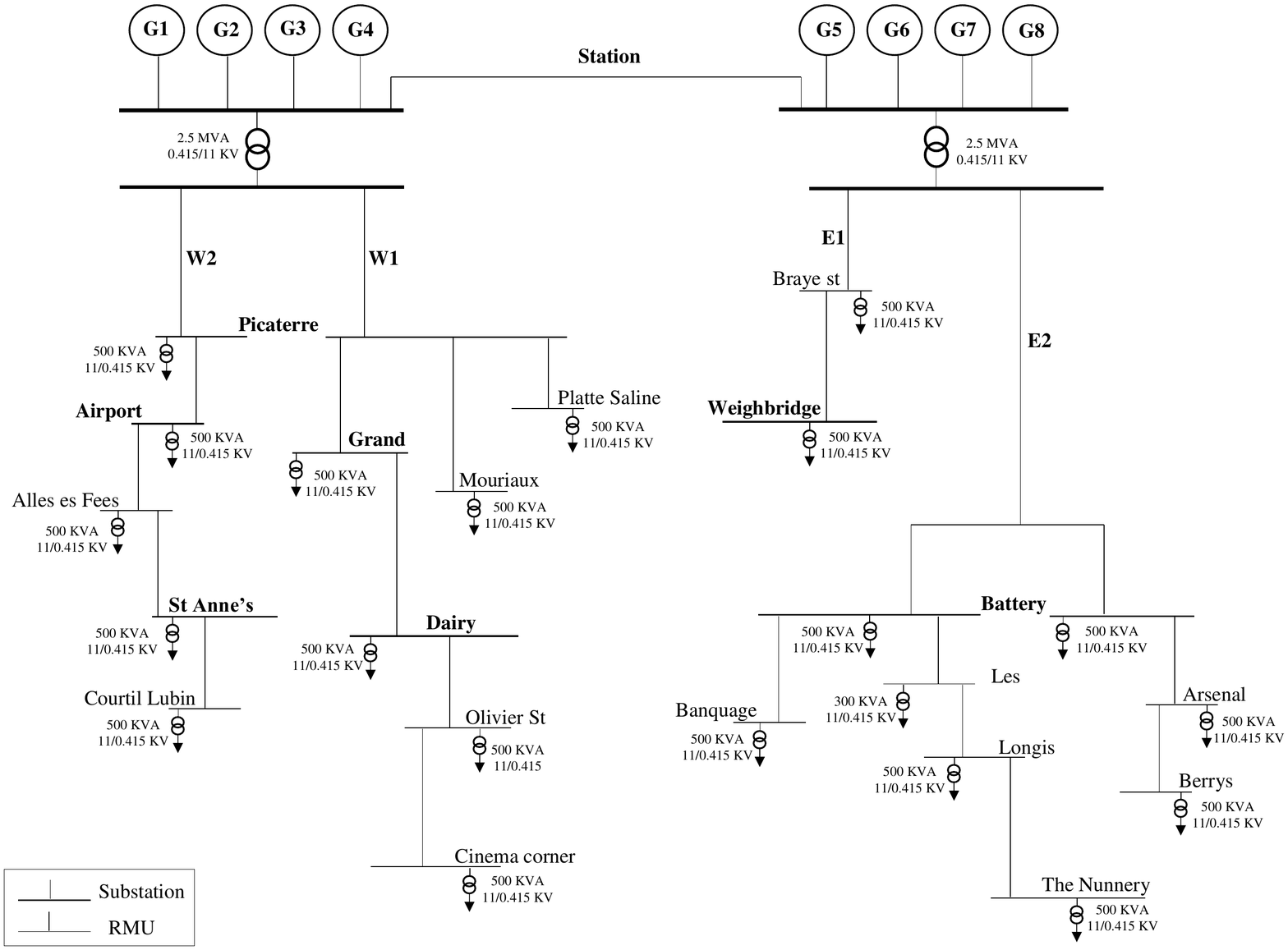}}
	\caption{The AEL network one-line diagram.}
	\label{AELold}
	\vspace{-0.5cm}
\end{figure} 
\begin{figure*}[!ht]
     \centering
     \begin{subfigure}[b]{0.3\textwidth}
         \centering
         \includegraphics[width=\textwidth]{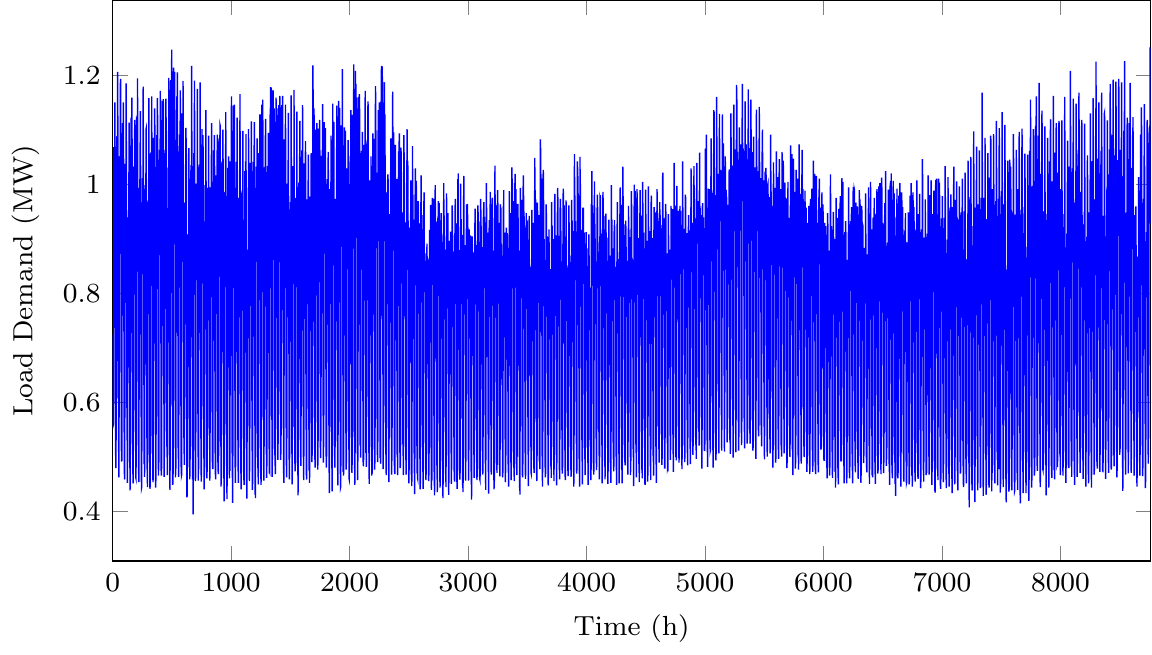}
     \end{subfigure}
     \hfill
     \begin{subfigure}[b]{0.3\textwidth}
         \centering
         \includegraphics[width=\textwidth]{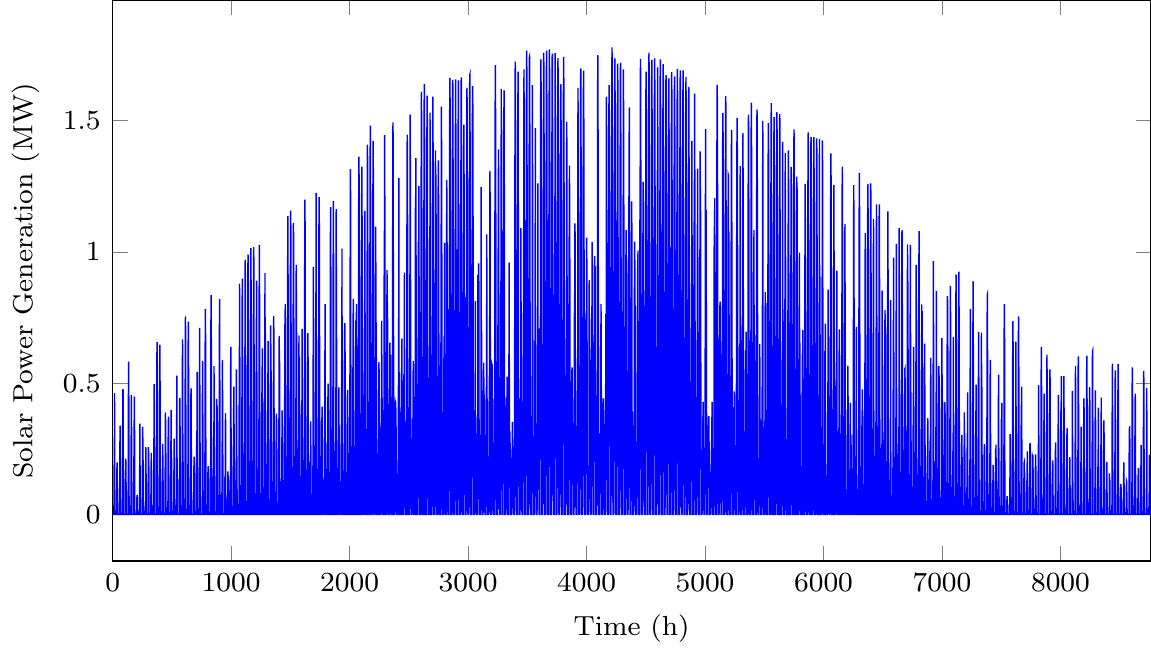}
     \end{subfigure}
     \hfill
     \begin{subfigure}[b]{0.3\textwidth}
         \centering
         \includegraphics[width=\textwidth]{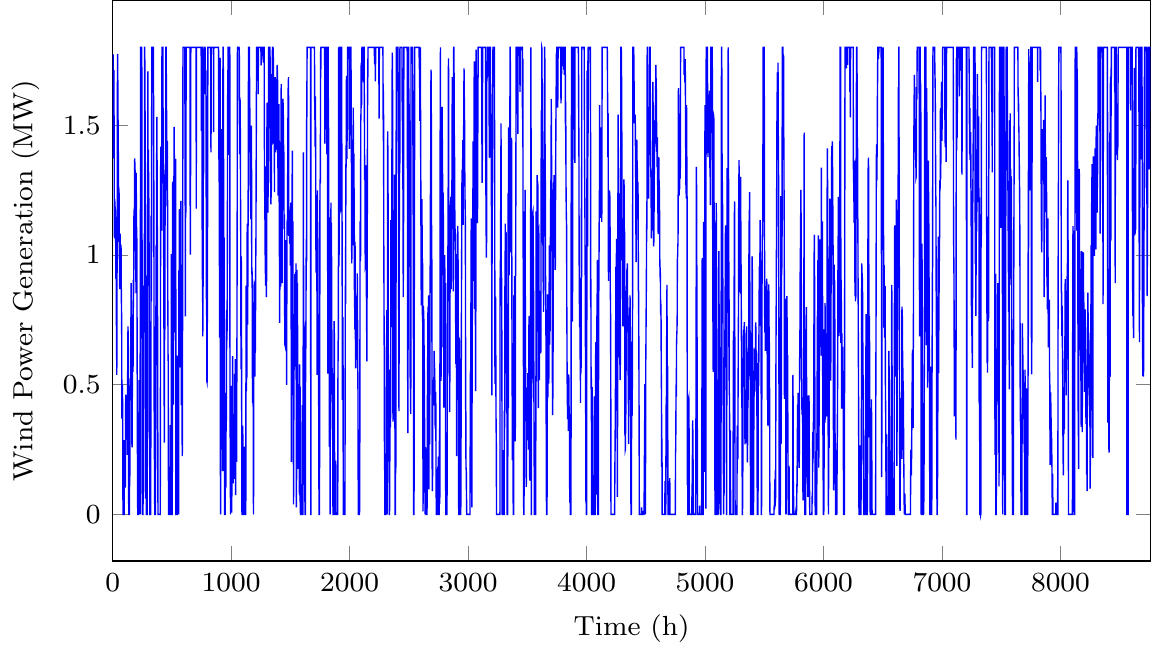}
     \end{subfigure}
        \caption{Yearly profiles of load demands and solar/wind power generations on Alderney island in 2013.}
        \label{fig:allprofiles}
\end{figure*}
\section{Case Studies}\label{section4}
\vspace{-.25cm}
\subsection{Input Data} 
\vspace{-.25cm}
In this section, the SMIP model described above is solved using PyEPLAN~\cite{shahab_dehghan_2020_3894705} to propose a low-carbon MG design for the Alderney electricity network. The AEL 11~kV primary network consists of four radial feeders as depicted in Fig.~\ref{AELold}. Electric power is generated solely at the power station by the $8 \times 450$~kW diesel units. The power station is connected to the 11~kV primary distribution network via two 2500~kVA transformers and the 11~kV primary distribution network is connected to the 415~V secondary distribution network by 500~kVA transformers at different substations and locations. The AEL distribution network comprises mainly three types of underground copper core cables (16 mm$^2$ PILC cables, 25 mm$^2$ PILC cables, and 70 mm$^2$ XLPE cables). There are a variety of other types and sizes of cable in certain locations. For example, newer additions to the high-voltage side are usually 70, 90 or 150 mm$^2$ cables.

Furthermore, battery, solar, and wind units are considered as investment candidates while investment costs of different technologies are taken from \url{https://atb.nrel.gov} and depicted in Table~\ref{investmentcosts}. Also, it is assumed that the interest rate (i.e., $i$) is equal to $0.053$, while the life time (i.e., $y$) of battery, solar, and wind units is equal 15, 30, and 30, respectively. Accordingly, the capital recovery factor (i.e., ${CRF} = \frac{i\cdot (1+i)^y}{(1+i)^y-1}$) for battery, solar, and wind units is equal to 0.098, 0.067, and 0.067, respectively, and consequently, the annualised investment costs can be calculated as depicted in Table~\ref{investmentcosts}. Also, it is assumed that operational costs of battery, solar, and wind units are equal to zero while the operational cost of diesel units is equal to 196.2 £/MWh~\cite{AEL_Website} on Alderney island at the time of writing, but fluctuates with market price on the date of loading at the refinery. The penalty cost of curtailing load demand is set to 1962 £/MWh.

The $k$-means clustering technique is used to obtain representative days using the yearly profiles of load demands and RES power generations on Alderney island in 2013. The peak load is equal to 1.252 MW. In addition, the solar irradiation and wind speed on Aldenery island in 2013 are taken from \cite{meteoblue}. In this paper, it is assumed that the efficiency of candidate solar panels/modules in solar farm is equal to 10\%~\cite{denholm2007regional} and the cut-in speed, rated speed, and cut-out speed of candidate wind turbines (i.e., Vestas V90 1.8 MW) are equal to 4~m/s, 12~m/s, and 25~m/s, respectively. In addition, the hub height of each wind turbine is equal to 80~m. Given a 1.8~MW solar farm with a 2-hectare land used to construct this power plant and a 1.8~MW wind farm, the yearly profiles of load demands, solar power generations, and wind power generations in 2013 are depicted in Fig.~\ref{fig:allprofiles}. 

\begin{table}[!t]
\centering
\caption{Investment costs of different technologies}
\label{investmentcosts}
\resizebox{\linewidth}{!}{%
\begin{tabular}{|c|c|c|c|} 
\hline
{Technology}                        & {Battery ({\textrm B})} & {Solar ({\textrm S})} & {Wind ({\textrm W})}  \\ 
\hline
{Investment Cost (M£/MW)}           & 0.98                   & 0.84                & 1.21                \\ 
\hline
{Annualised Investment Cost (£/MW)} & 96040                  & 56280               & 81070               \\
\hline
\end{tabular}
}
\end{table}
\begin{table}[!t]
\vspace{-0.25cm}
\centering
\caption{RES capacity factor on Alderney island}
\label{capacityfactor}
\resizebox{1\linewidth}{!}{%
\begin{tabular}{|c|c|c|} 
\hline
{Technology} & {Built Capacity (MW)} & {Capacity Factor (\%)}  \\ 
\hline
{\textrm Solar (S)}               & 1.8                              & 16.27                          \\ 
\hline
{\textrm Wind (W)}                & 1.8                              & 54.39                          \\
\hline
\end{tabular}
}
\end{table}

The capacity factors\footnote{The capacity factor represents the ratio of the electrical energy generated by a specific technology to the electrical energy, which could have been generated at rated capacity continuously during a one-year period (or other specific periods).} (CFs) of both solar and wind farms are presented in Table~\ref{capacityfactor}. Accordingly, the CF of wind technology is significantly higher than the CF of solar technology while the land needed by wind turbines to create a 1.8 MW wind farm is significantly less of the land needed by solar panels/modules to create a 1.8 MW solar farm (i.e., approximately 2 hectares). Additionally, according to Table~\ref{investmentcosts}, battery units have the highest annualised investment costs while solar units have the lowest annualised investment costs. Therefore, it is necessary to use the proposed planning tool to obtain the optimal technology mix for creating a sustainable MG on Alderney island under different circumstances.

\begin{figure*}[!ht]
     \centering
     \begin{subfigure}[b]{0.3\textwidth}
         \centering
         \includegraphics[width=\textwidth]{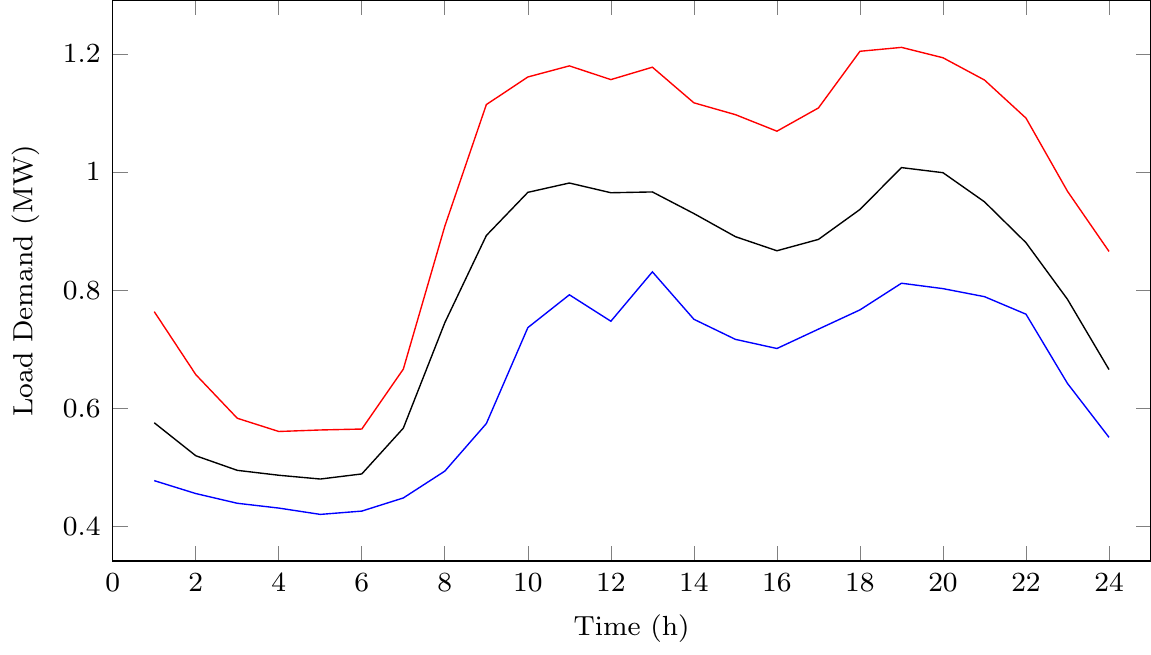}
     \end{subfigure}
     \hfill
     \begin{subfigure}[b]{0.3\textwidth}
         \centering
         \includegraphics[width=\textwidth]{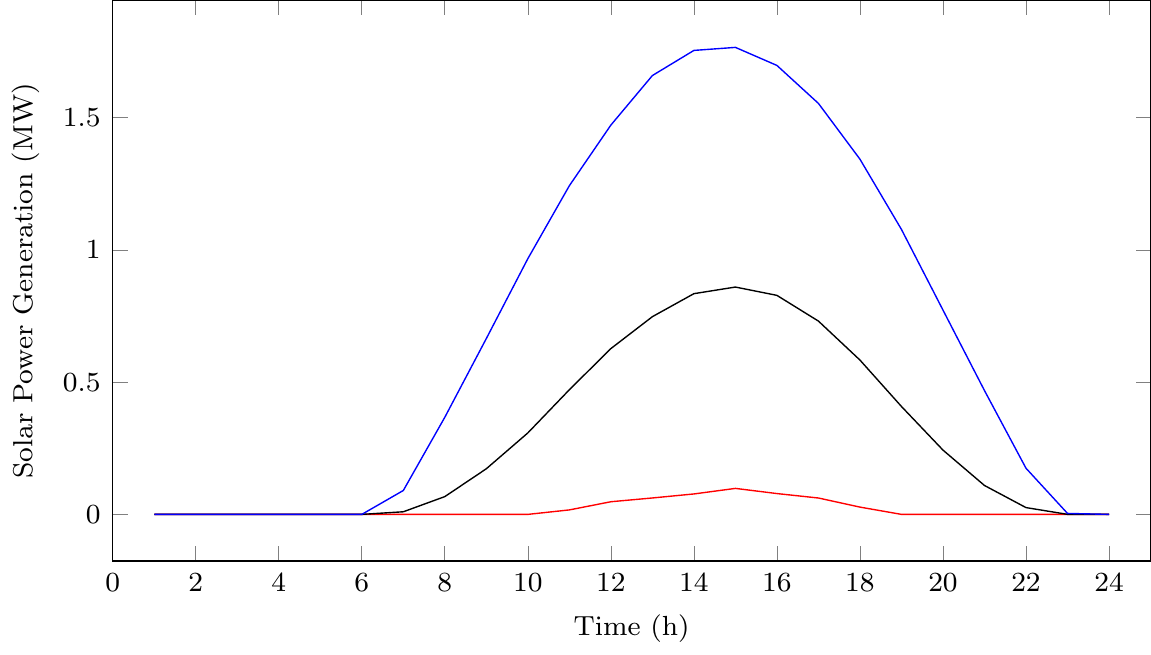}
     \end{subfigure}
     \hfill
     \begin{subfigure}[b]{0.3\textwidth}
         \centering
         \includegraphics[width=\textwidth]{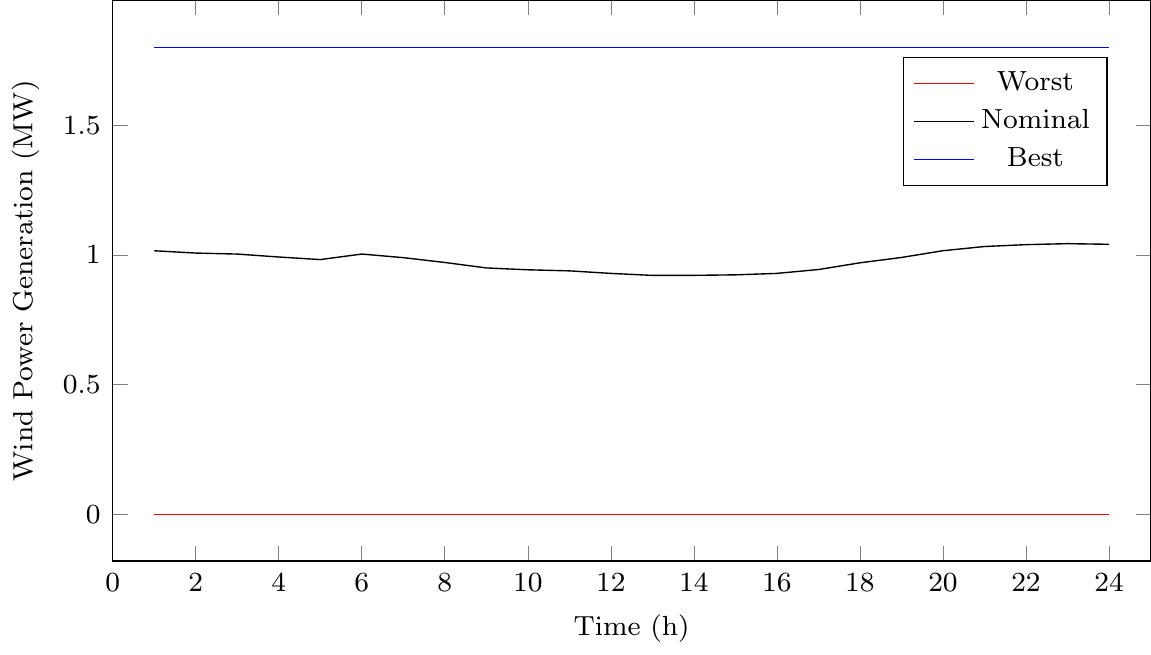}
     \end{subfigure}
        \caption{The best, nominal, and worst representative days for load demands and solar/wind power generations on Alderney island.}
        \label{fig:1repday}
        \vspace{-0.25cm}
\end{figure*}
\begin{table*}[!ht]
\centering
\caption{Optimal investment plans for different cases under best, nominal, and worst representative days}
\label{allcasest}
\begin{tabular}{|p{5cm}|P{1.76cm}|P{1.76cm}|P{1.76cm}|P{1.76cm}|P{1.76cm}|P{1.76cm}|} 
\hline
{Case Number}                & {C1}                           & {C2}                           & {C3}                           & {C4}                                                        & {C5}                                                        & {C6}                            \\ 
\hline
{Best Representative Day}    & $1\times \textrm{W}$ & $1\times \textrm{W}$ & $5\times \textrm{S}$ & $1\times \textrm{B}$,$2\times \textrm{S}$ & $1\times \textrm{W}$                              & $1\times \textrm{W}$  \\ 
\hline
{Nominal Representative Day} & $2\times \textrm{W}$ & $2\times \textrm{W}$ & $9\times \textrm{S}$ & $1\times \textrm{B}$,$5\times \textrm{S}$ & $1\times \textrm{S}$,$1\times \textrm{W}$ & $1\times \textrm{W}$  \\ 
\hline
{Worst Representative Day}   & Infeasible                            & Infeasible                            & $10\times \textrm{S}$     & $10\times \textrm{S}$                             & $10\times \textrm{S}$                             & AEL MG                    \\
\hline
\end{tabular}
\vspace{-0.25cm}
\end{table*}
\begin{figure*}[!ht]
     \centering
     \begin{subfigure}[b]{0.3\textwidth}
         \centering
         \includegraphics[width=\textwidth]{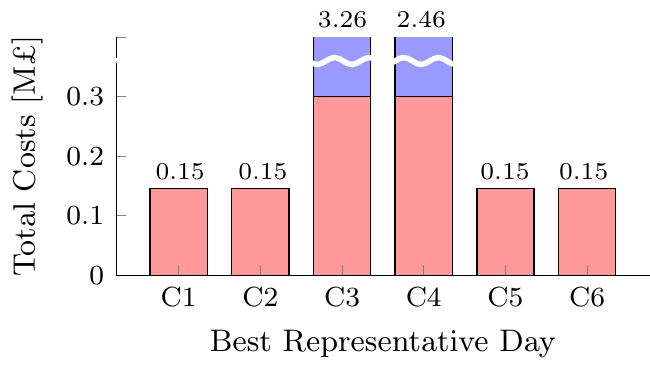}
         \caption{Risk-Seeker}
         \label{risk-seeker}
     \end{subfigure}
     \hfill
     \begin{subfigure}[b]{0.3\textwidth}
         \centering
         \includegraphics[width=\textwidth]{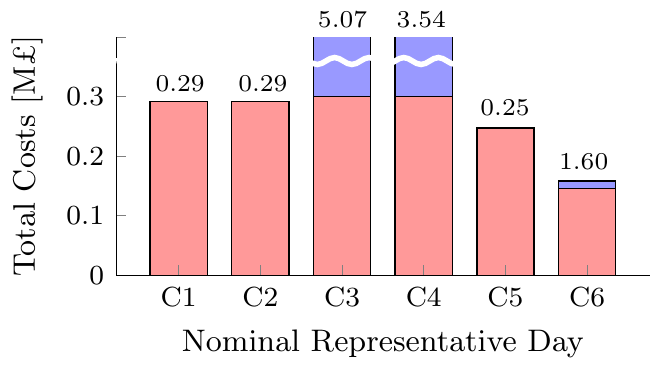}
         \caption{Risk-Neutral}
         \label{risk-neutral}
     \end{subfigure}
     \hfill
     \begin{subfigure}[b]{0.3\textwidth}
         \centering
         \includegraphics[width=\textwidth]{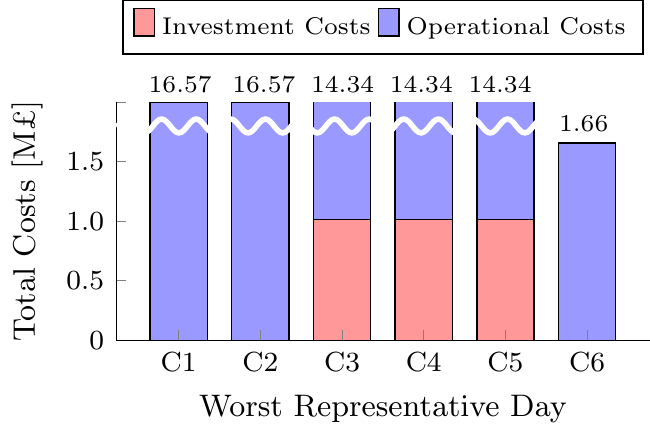}
         \caption{Risk-Averse}
         \label{risk-averse}
     \end{subfigure}
        \caption{Total investment and operational costs for different cases under best, nominal, and worst representative days.}
        \label{fig:allcases}
        \vspace{-0.25cm}
\end{figure*}
\vspace{-0.25cm}
\subsection{Investment Plan Under Best, Nominal, and Worst Scenarios} 
\vspace{-0.25cm}
In this study, one best, nominal, and worst representative day are constructed using the yearly profiles of load demands and solar/wind power generations on Alderney island in 2013, as illustrated in Fig.\ref{fig:allprofiles}. Also, different investment alternatives are considered at the current location of the AEL power plant, including: 
\noindent
\textbf{Case 1 (C1):} Only $10 \times 1.8$-MW wind units are considered as investment candidates.

\noindent
\textbf{Case 2 (C2):} Both $10 \times 1.8$-MW battery units and $10 \times 1.8$-MW wind units are considered as investment candidates. 

\noindent 
\textbf{Case 3 (C3):} Only $10 \times 1.8$-MW solar units are considered as investment candidates. 

\noindent 
\textbf{Case 4 (C4):} Both $10\times  1.8$-MW battery units and $10 \times 1.8$-MW solar units are considered as investment candidates. 

\noindent 
\textbf{Case 5 (C5):} All $10 \times 1.8$-MW battery units, $10\times 1.8$-MW solar units, and $10\times 1.8$-MW wind units are considered as investment candidates. 

\noindent 
\textbf{Case 6 (C6):} In addition to the current AEL diesel units, all options in C5 are considered as investment candidates in C6.

The best, nominal, and worst representative days are illustrated in Fig.~\ref{fig:1repday} wherein solar/wind power generations are provided for each unit. The optimal investment plans for all cases under the best, nominal, and worst representative days are presented in Table~\ref{allcasest}. Moreover, the total investment and operational costs are depicted in Fig.~\ref{fig:allcases}. For all cases C1-C6, the total costs under the best representative day have the lowest value while the total costs under the worst representative day have the highest value. For instance, the total costs for the best, nominal, and worst representative days are equal to $0.15$ M£ in Fig.~\ref{risk-seeker}, $0.29$ M£ in Fig.~\ref{risk-neutral}, and $16.57$ M£ in Fig.~\ref{risk-averse}, respectively. It is noteworthy to mention that the best representative day for wind power generation corresponds to the maximum capacity of each candidate wind unit while the worst representative day for wind power generation corresponds to no power generation. Accordingly, C1 and C2 under the worst representative day result in {\it infeasible} solutions, as illustrated in Table~\ref{allcasest}, and their total costs in  Fig.~\ref{risk-averse} (i.e., $16.57$ M£) only correspond to the penalty cost of load demand curtailment during the entire planning horizon. However, C1, C2, C5, and C6 under the best representative day result in identical optimal investment plans, only constructing a $1.8$ MW wind unit and obviating the need to operate the current AEL diesel units. Furthermore, C6 provides not only the lowest total costs, similar to C1, C2, C5, and C5, under the best representative day, but also the lowest total costs under the nominal and worst representative days. However, C6 under the worst representative day only rely on the current AEL MG without constructing any battery, solar, or wind units. The main reason is that creating a sustainable MG on Alderney based on only one worst representative day results in an over-conservative investment plan.                   
\begin{table*}[!ht]
\centering
\caption{Investment plans for different number of best, nominal, and worst representative days for Case~C6}
\label{allrdayst}
\begin{tabular}{|p{5cm}|P{2.19cm}|P{2.19cm}|P{2.19cm}|P{2.19cm}|P{2.19cm}|} 
\hline
{Case Number}                & {R1}                    & {R5}                                                                                    & {R10}                                       & {R50}                                        & {R100}                                        \\ 
\hline
{Best Representative Day}    & $1\times \textrm{W}$  & $1\times \textrm{W}$                                                                  & $1\times \textrm{W}$                      & $1\times \textrm{W}$                       & $1\times \textrm{S}$,$1\times \textrm{W}$   \\ 
\hline
{Nominal Representative Day} & $1\times \textrm{W}$  & \begin{tabular}[c]{@{}l@{}}$1\times \textrm{S}$,$1\times \textrm{W}$ \end{tabular} & $1\times\textrm{S}$,$1\times\textrm{W}$ & $1\times \textrm{S}$,$1\times\textrm{W}$  & $1\times \textrm{S}$,$1\times \textrm{W}$   \\ 
\hline
{Worst Representative Day}   & AEL MG            & $2\times \textrm{W}$                                                                  & $1\times \textrm{S}$,$1\times \textrm{W}$ & $1\times \textrm{S}$,$1\times \textrm{W}$  & $1\times \textrm{S}$,$1\times \textrm{W}$   \\
\hline
\end{tabular}
\end{table*}
\begin{figure*}[!ht]
     \centering
     \begin{subfigure}[b]{0.3\textwidth}
         \centering
         \includegraphics[width=\textwidth]{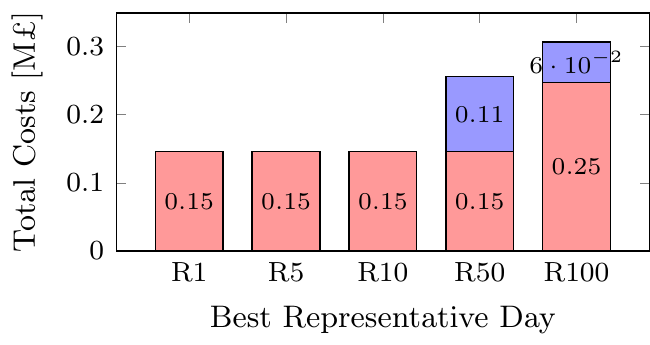}
         \caption{Risk-Seeker}
         \label{risk-seeker2}
     \end{subfigure}
     \hfill
     \begin{subfigure}[b]{0.3\textwidth}
         \centering
         \includegraphics[width=\textwidth]{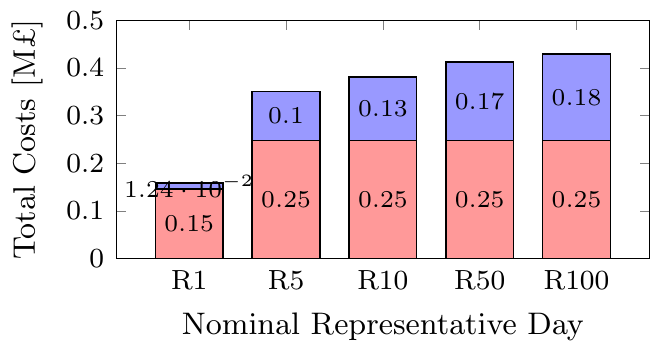}
        \caption{Risk-Neutral}
         \label{risk-neutral2}
     \end{subfigure}
     \hfill
     \begin{subfigure}[b]{0.3\textwidth}
         \centering
         \includegraphics[width=\textwidth]{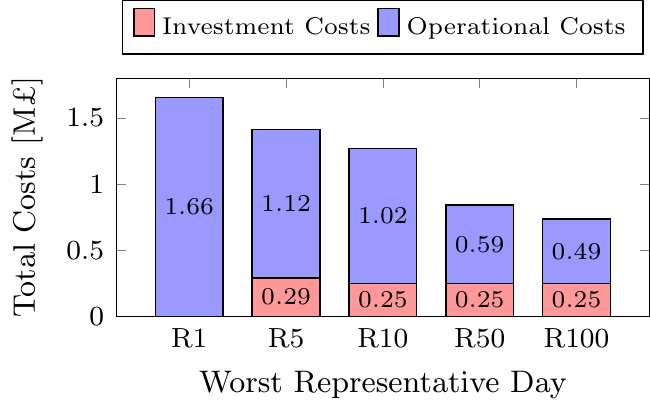}
        \caption{Risk-Averse}
         \label{risk-averse2}
     \end{subfigure}
        \caption{Total investment and operational costs for different number of best, nominal, and worst representative days.}
        \label{fig:alldays}
        \vspace{-0.25cm}
\end{figure*}

\subsection{Investment Plan for Different Number of Representative Days} 
\vspace{-0.25cm}
To enhance the accuracy of the proposed solution, different number of best, nominal, and worst representative days can be considered for C6, including $1$ (R1), $5$ (R5), $10$ (R10), $50$ (R50), and $100$ (R100). The optimal investment plan for C6 for each choice are presented in Table~\ref{allrdayst} and their total investment and operational costs are depicted in  Fig.~\ref{fig:alldays}. Increasing the number of representative days increases the total costs under the best representative day (Fig.~\ref{risk-seeker2}) and the nominal representative day (Fig.~\ref{risk-neutral2}), while decreases the total cost under the worst representative day (Fig.~\ref{risk-averse2}). Additionally, the investment plans are identical under the best, nominal, and worst representative days in R100 (constructing one $1.8$ MW solar and one $1.8$ MW wind unit in addition to the current AEL MG). It is worthwhile to mention that the optimal investment plan remains unchanged after $5$ representative days under the nominal condition, while it remains unchanged after $100$ (resp. $10$) representative days under the best (resp. worst) conditions, as shown in Fig.~\ref{fig:alldays}. Finally, it can be concluded that $5$ nominal representative days can appropriately characterise the uncertain profiles of load demand and RES generation on Alderney island with reasonable computational complexity.   

\section{Conclusion}\label{section5}
This paper presents a two-stage stochastic model for creating a sustainable MG on Alderney island under the uncertainty of load demands and RES power generations. Also, the $k$-means clustering technique is used to characterise the yearly profiles of load demands and RES power generations through a sufficient number of best, nominal, and worst representative days. The proposed MG planning model is implemented in the open-source tool PyEPLAN. Simulation results demonstrate that the best low-carbon investment plan pertains to a hybrid MG including both solar and wind power in addition to current AEL diesel units.   

\section{References}\label{refrences}
\bibliographystyle{iet.bst}
\balance
{\bibliography{refs}}
\end{document}